\documentclass[12pt]{article}%
\usepackage{graphicx}
\usepackage{amsmath}%
\usepackage{amsfonts}%
\usepackage{amssymb}
\newtheorem{theorem}{Theorem}[section]

\newtheorem{definition}[theorem]{Definition}

\newtheorem{proposition}[theorem]{Proposition}
\newtheorem{remark}[theorem]{Remark}

\newenvironment{proof}[1][Proof]{\textbf{#1.} }{\ \rule{0.5em}{0.5em}}

\begin{document}

\author{Gadde A. Swarup\\Department of Mathematics and Statistics\\University
of Melbourne\\Victoria 3010, Australia.\\email:gadde@ms.unimelb.edu.au }
\title{Delzant's variation on Scott Complexity}
\maketitle
\begin{center}
\textbf{To Peter Scott, for his 60th birthday}
\end{center}

\begin{abstract}
We give an exposition of Delzant's ideas extending the notion of Scott 
complexity of finitely generated groups
to surjective homomorphisms from finitely presented groups.
\end{abstract}

\newpage

\newpage

\textbf{{\Large Introduction}}

It is an old question of W. Jaco \cite{Jaco} whether every finitely generated 
indecomposable group
has a finitely presented indecomposable cover: given a finitely generated
indecomposable group $H $
is there a surjective homomorphism $\phi:G\rightarrow H$ from a
finitely presented group $G$ such that
for any factorization of $\phi=\psi \alpha$, 
$G\stackrel{\alpha}{\rightarrow} G'\stackrel{\psi}{\rightarrow} H$ 
with $\alpha$ surjective,
we have $G'$ is also indecomposable. Jaco originally raised the 
question in connection with the coherence of 3-manifold groups.
Peter Scott \cite{Scott} (and independently Peter Shalen) proved the coherence 
of $3$-manifold groups but bypassed the above question.
In his proof, Scott used a notion of complexity of finitely 
generated groups which is sometimes called Scott complexity.
Thomas Delzant extended the notion of complexity to surjective homomorphisms 
of finitely presented groups to
finitely generated groups and answered the above question in the 
affirmative. This leads to a quick proof
of the coherence of $3$-manifold groups as well as a proof of the 
acylindrical accessibility theorem of Sela.
Delzant knew these arguments for several years and seems to have other applications in mind.
Hopefully, he will write up a more complete exposition of his ideas. In view of the interest
shown by various people who worked on this problem, we 
give an exposition of some of his arguments.

\medskip
\section{Complexity}

We call a group $G$ {\em indecomposable} if $G$ is not a free product 
and is not isomorphic to ${\mathbb Z}$. Some authors use the term
`freely
indecomposable'. Let $G$ be a finitely generated group and let 
$$G=G_1*\cdots *G_m*G_{m+1}*\cdots *G_{m+n}$$ be
a free product  decomposition of $G$ with $G_i$ indecomposable for 
$i\leq m$ and $G_i$ isomorphic to ${\mathbb Z}$
for $i>m$. The factors $G_i,i\leq m$ are called the {\em indecomposable 
factors} of $G$. These are
unique up to isomorphism and there are only finitely many such up to 
conjugacy in $G$. We call decompositions
of $G$ of the above type standard decompositions of $G$. The ordered 
pair $(m+n,n)$ is called the complexity or Scott
complexity of $G$ and is denoted by $c(G)$. There is a partial order 
on the complexities given by lexicographic
order. Scott used Stallings technique of binding ties \cite{Stallings} to prove:
\begin{theorem}
Let $\phi :G\rightarrow H$ be a surjective homomorphism of finitely 
generated groups such that
$\phi $ restricted to each indecomposable factor of $G$ is injective. 
Then $c(G)\geq c(H)$ and
c(G)=c(H) if and only if $\phi$ is an isomorphism.
\end{theorem}
This theorem is proved by first taking a standard decomposition 
$$H=H_1*\cdots *H_m*H_{m+1}*\cdots *H_{m+n}$$
where $C(H)=(m+n,n)$ and then obtaining a decomposition 
$$G=G_1*\cdots *G_m*G_{m+1}*\cdots *G_{m+n}$$ with
$\phi (G_i)\subseteq H_i$. This is achieved by the method of binding 
ties. The hypothesis that
$\phi$ is injective on the indecomposable factors of $G$ implies that the
$G_j,j\geq m$, are free.
Thus $c(G)\geq c(H)$. If $c(G)=c(H)$, clearly $\phi$ is an isomorphism.

The above argument easily extends to the case when $\phi (K)$ is 
indecomposable for each indecomposable factor
of $G$. More generally:
\begin{theorem}
Let $\phi :G\rightarrow H$ be a surjective homomorphism of finitely 
generated groups such that
$\phi (K)$ can be conjugated into an indecomposable factor of $H$
for each indecomposable factor $K$ of $G$. Then $c(G)\geq c(H)$. Suppose
that $c(G)=c(H)=(m+n,n)$ and $$H=H_1*\cdots *H_m*H_{m+1}*\cdots *H_{m+n}$$
is a  standard decomposition of $H$. Then
there is a standard decomposition $$G=G_1*\cdots
*G_m*G_{m+1}*\cdots *G_{m+n}$$  such that $\phi (G_i)=H_i$.
\end{theorem}
Thus, when $c(G)=c(H)$, and $\phi (K)$ are indecomposable
for each indecomposable factor $K$ of $G$, the 
standard decompositions of $H$ can
be imitated by standard decompositions of $G$ which respect $\phi$. 
It is also easy to see that in this case,
the standard decompositions of $G$ can be pushed forward.
\begin{remark}
Suppose that $\phi:G\rightarrow H$ is a surjective homomorphism of 
finitely generated groups such that $\phi (K)$ is indecomposable for
each indecomposable factor $K$ of 
$G$ and further assume
that $c(G)=c(H)=(m+n,n)$. If $$G=G_1*\cdots *G_m*G_{m+1}*\cdots *G_{m+n}$$
is a  standard decomposition of $G$, then there is a
standard decomposition $$H=H_1*\cdots *H_m*H_{m+1}*\cdots *H_{m+n}$$ such
that 
$\phi (G_i)=H_i$.
\end{remark}
To see this start with a standard decomposition 
$$H=H'_1*\cdots *H'_m*H'_{m+1}*\cdots *H'_{m+n}.$$Consider $G'=*\phi
(G_i)$. We have surjective homomorphisms:
  $$G\stackrel{\alpha}{\rightarrow} G'\stackrel{\psi}{\rightarrow} H$$
with $\phi=\psi \alpha$. We also have $c(G')=(m+n,n)$.
By the previous theorem, there is a standard decomposition 
$G'=G'_1*\cdots *G'_m*G'_{m+1}*\cdots *G'_{m+n}$
with $\alpha (G'_i)=H'_i$. By construction $\alpha $ restricted to 
the indecomposable
factors of $G'$ is injective and thus $\alpha $ is
an isomorphism. Thus we can take $H_i=\alpha \psi (G_i)$.

We now give Delzant's extension of the notion of complexity.
\begin{definition}
Let $\phi :G\rightarrow H$ be a surjective homomorphism of 
a finitely presented group G. Consider
factorizations of $\phi=\psi \alpha$, 
$G\stackrel{\alpha}{\rightarrow} G'\stackrel{\psi}{\rightarrow} H$
with $\alpha $ surjective and $G'$ finitely presented. The complexity 
of $\alpha$ is by definition the
supremum of the complexities $c(G')$ as $G'$ varies over finitely 
presented groups.
\end{definition}
We also mention a folklore result that is used below:
\begin{proposition}
Let $\phi :G\rightarrow H$ be an epimorphism of a finitely presented group $G$
onto a finitely generated group $H$. Suppose that $H$ is a non-trivial free product $H_1*H_2$.
Then there is a factorization $G\stackrel{\alpha}{\rightarrow} G'\stackrel{\psi}{\rightarrow} H$
with $G'$ finitely presented, $\alpha$ surjective, $G'=G_1*G_2$ and $\psi
(G_i)=H_i,i=1,2$.
\end{proposition}
This proposition is easily proved by Stallings' method of binding ties.

In the above discussion, we considered several times homomorphisms $\phi :G\rightarrow H$ which do not
have factorization of the form $\phi=\alpha \psi$, where $\psi$ is a surjective homomorphism from
$G$ to $G'$ with $G'$ either a non-trivial free product or infinite cyclic. It is natural to call
such homomorphisms $\phi$  {\em essential }.
\section{Main Results}
\begin{theorem}\label{delzant1}
Let $H$ be a finitely generated indecomposable group. Then there is a
finitely presented indecomposable group $G$
and a surjective homomorphism $\phi :G\rightarrow H$ such that for 
any factorization: $\phi=\psi \alpha$, 
$$G\stackrel{\alpha}{\rightarrow} 
G'\stackrel{\psi}{\rightarrow} H$$
with $\alpha $ surjective, then $G'$ is also 
indecomposable.
\end{theorem}
Such a $\phi :G\rightarrow H$ is called a finitely presented indecomposable cover of $H$.
In the terminology introduced at the end of the previous section, 
Theorem \ref{delzant1} reads:
\begin{theorem}\label{delzant2}
If $H$ is a finitely generated indecomposable group, then $H$ admits
essential epimorphisms from finitely presented groups.
\end{theorem}
In view of Proposition 1.5, Theorem \ref{delzant2} easily follows from:

\begin{theorem}
Let  $\phi _i:G_i \rightarrow H$,
$G_i\stackrel{\alpha _i}{\rightarrow} G_{i+1}\stackrel{\phi 
_{i+1}}{\rightarrow} H$
be such that $\phi _i =\phi _{i+1} \alpha _i$ and all maps are 
surjective homomorphisms with
$G_i$ finitely presented. Suppose further that $G$ is the direct limit 
of the $G_i$, that is, if $\phi _1 (g)=1$,
then there is an $n$ such that $\alpha _n \alpha _{n-1} ...\alpha _1 
(g)=1$. Then there is an integer $K$ such
that $G_i$ is indecomposable for $i\geq K$.
\end{theorem}
We now present the proof of the second theorem. 

\noindent\begin{proof}
Clearly 
$c(\phi _i)\geq c(\phi _{i+1})$. By going to a subsequence
if necessary, we may assume that $c(\phi _i)$ are all equal to $(m+n,n)$. If
$(m+n,n)=(1,0)$, there is nothing to prove. Otherwise, we will arrive 
at a contradiction.
Consider standard decompositions of $G_i$.
$$G_i=G^i_1*\cdots *G^i_m*G^i_{m+1}*\cdots *G^i_{m+n}.$$ 
We want to arrange 
these so that $\alpha _i (G^i _j)\subseteq G^{i+1}_j$.
We claim that $\alpha _i (G^i _j)$ is indecomposable. Firstly
$\alpha _i (G^i _j)$ cannot be
isomorphic to ${\mathbb Z}$ since $c(G_i)=c(g_{i+1})$. If $\alpha _i (G^i
_j)$  is a free product, then $\alpha _i$
factors through a finitely presented group with at least $(n+m+1)$ 
factors by Proposition 1.5. Hence
$\alpha _i(G^i _j)$ are indecomposable for all $i$ and $j\leq m$. 
Hence, by Remark 1.3, we can arrange
the standard decompositions of $G_i$ so that $\alpha _i (G^i _j)= 
G^{i+1}_j$. To complete the proof of
Theorem 2.2, we observe that since $G$ is a direct limit of $G_i$, 
$G$ is a free product
of the direct limits $G^j$ of $G_j^i$. Since $G$ is indecomposable, 
all but one of $G^j$ must be trivial.
But none of $G^j$ can be trivial since all $G^i _j$ are finitely 
presented and triviality of $G^j$ implies that
the complexity of some $G_i$ for large $i$ is smaller than 
$(m+n,n)$. Clearly $(m+n,n)\neq (0,1)$. This completes
the proof of Theorem \ref{delzant2}.
\end{proof}

\section{Applications}
We sketch quick proofs of two applications.
The first is the Scott-Shalen theorem \cite{Scott}:
\begin{theorem}
If the fundamental group of a 3-manifold is finitely generated, then 
it is finitely presented.
\end{theorem}
The argument goes as follows. Let $H=\pi _1 (M)$ be the fundamental 
group of a $3$-manifold $M$ and suppose that $H$ is finitely generated.
We may assume that $H$ is indecomposable.
Let $\phi :G\rightarrow H$ be a finitely presented indecomposable cover of $H$.
We represent $G$ as the fundamental group of a finite simplicial
complex $K$ and construct a piecewise linear map $f$ which induces $\phi$.
Let $N_1$ be a regular neighbourhood of $f(K)$ and $G_1$ be the image
of $\pi _1 (K)$ in $\pi _1 (N_1)$. By Theorem 1.1, $G_1$ is indecomposable.
If the boundary $\partial N_1$ of $N_1$ is not incompressible in $M$, then there is a Dehn disc $D_1$
contained in either $N_1$ or in $C_1$, the closure of the complement of $N_1$ in $M$ such that $D_1$
intersects $\partial N_1$ in exactly $\partial D_1$. If $D_1$ is contained in $C_1$, we add thickened $D_1$ to
$N_1$ to obtain $N_2$ and call $G_2$, the image of $G$ in $\pi _1 (N_2)$. If $D_1$ is in $N_1$, we split $N_1$ along $D_1$.
Then one of the pieces (one may be empty if $D_1$ is non-separating), say $N_2$ contains $G_1$ up to
conjugacy. Call this $G_2$. So, we can homotope $f$ so that the image of $G$ is $G_2$ under the induced map in the fundamental groups.
After a finite number of steps, we find $N_k$ which is incompressible in $M$. Hence $\pi _1 (N_k)$ maps injectively to $\pi _1 (M)$. 
Since $\pi _1 (N_k)$ contains $G_k$ which maps onto $H=\pi _1 (M)$, we see that $\pi _1 (N_k)$ is isomorphic to $\pi _1(M)$. Since $N_k$ is compact,
we see that $\pi _(M)$ is finitely presented.
 
The second application is Sela's acylindrical accessibility theorem (see  \cite{Sela} and \cite{Weid}).
Delzant has given an elementary proof of a more general result in the finitely presented case \cite{Delzant}.
A simplicial action of a group $H$ on a simplicial tree $T$ is said 
to be $k$-acylindrical, if the stabilizers
of segments of length $ k$ are trivial.
\begin{theorem}[Sela]\label{sela1}
Let $H$ be a finitely generated indecomposable group and $k$ a 
positive integer. Then there
is a number $n(k,H)$ such that for any $k$-acylindrical minimal 
action of $H$ on a simplicial tree $T$,
the number of vertices of $T/H$ is bounded by $n(k,H)$.
\end{theorem}

We recall Delzant's generalization in the finitely presented case.
\begin{definition}
Let ${\mathcal C}$ be a family of subgroups of  group $G$ which is closed 
under conjugation and taking subgroups. We say that a $G$-tree $T$ is $(k,{\mathcal C})$-acylindrical
if the stabilizers of segments of length $k$ are in ${\mathcal C}$.
\end{definition}
Delzant shows:
\begin{theorem}[Delzant \cite{Delzant0}] \label{delzant3}
Suppose $G$ is a finitely presented group and ${\mathcal C}$ is a family of subgroups of $G$ which is closed
under conjugation and taking subgroups. Moreover, suppose that $G$ does not split
over an element of ${\mathcal C}$. If $T$ is a minimal $(k,{\mathcal C})$-acylindrical tree, then there is a number $n(k,G)$
such that the number of vertices of $T/G$ is bounded by $n(k,G)$.
\end{theorem}
The number $n(k,G)$ is defined in terms of the triangular presentations of $G$.
In Weidmann's proof \cite{Weid}, the bound  is defined in terms of minimal number generators of $G$ . Delzant's argument gives
a more general result:
\begin{theorem}
Suppose that $G$ is a finitely presented group and ${\mathcal C}$ a family
of subgroups of $G$ closed under conjugation and taking subgroups. For
any positive integer $k$, there is a positive integer $n(k,G)$ such that
the following holds. For any $(k,{\mathcal C})$-acylindrical, minimal
$G$-tree $T$, $G$ has a graph of groups decomposition with edge groups in
${\mathcal C}$ (the decomposition may be trivial) such that for any of
the vertex groups $G_v$ then for a minimal $G_v$-subtree $T_v$ of $T$, $
T_v/ G_v$ has at most
$n(k,G)$ vertices. In particular, if $G$ does not split over an element of ${\mathcal C}$, then $T/G$
has less than $n(k,G)$ vertices.
\end{theorem}

The proof seems to be inspired by Dunwoody's ideas from  \cite{Dunwoody}. To deduce Sela's acylindrical accessibility theorem
from the above theorem,
let $\phi :G\rightarrow H$ be a finitely presented indecomposable cover of $H$
and take ${\mathcal C}$ to be the family of subgroups of the kernel of $\phi $.
Let $\Gamma $ be the graph of groups given by Delzant's theorem
\ref{delzant3}. Since every element of ${\mathcal C}$ fixes $T$, not all
the vertex groups of $\Gamma $ can be in ${\mathcal C}$. Choose a vertex
group $G_v$ with $\phi (G_v)\neq 1$. We claim that $\phi (G_v)=H$, for
otherwise $\phi $ factors through a free product. This completes the
proof of Sela's theorem \ref{sela1}.

Delzant's results and arguments seem to give a procedure 
for deducing results for finitely generated indecomposable groups from similar results in the finitely presented case.
Delzant has used similar ideas in \cite{Delzant} to study conjugacy
classes of homomorphic images of a finitely presented group in hyperbolic
group.

\end{document}